\documentclass[twocolumn,twoside]{article}
\usepackage{epsfig,amsmath,amsfonts}
\usepackage{amsmath,amssymb,amsfonts,amsthm}
\usepackage{enumerate}
\usepackage{verbatim}
\usepackage{fancyhdr}
\usepackage{pxfonts}
\usepackage{graphicx}
\usepackage{framed}
\usepackage{multicol}
\usepackage{mathrsfs}
\usepackage{epsfig}
\usepackage{leftidx}
\usepackage{color}
\usepackage{color}
\usepackage{appendix}
\usepackage[authoryear]{natbib}
\usepackage{fancyhdr}

\theoremstyle{plain}

\theoremstyle{definition}

\theoremstyle{remark}

\def\sech{\mathop{\operator@font sech}\nolimits}

\newtheorem*{thm*}{Theorem}

\newtheorem*{def*}{Definition}
\theoremstyle{example}

\newcommand{\red}{}

\newcommand{\R}{\mathbb{R}}
\newcommand{\C}{\mathbb{C}}

\newcommand{\G}{\mathcal{G}}

\renewcommand{\epsilon}{\varepsilon}
\renewcommand{\phi}{\varphi}

\begin{document}


\title{\bf Stochastic Parametrization using Compressed Sensing: Application to the Lorenz-96 Atmospheric Model}
\author{A. MUKHERJEE$^1$, Y. AYDOGDU$^1$, T. RAVICHANDRAN$^1$ and N. SRI NAMACHCHIVAYA$^1$
\thanks{Corresponding author.\hfil\break e-mail: navam@uwaterlooo.ca}\\
\small{$^1$Department of Applied Mathematics, University of Waterloo, Waterloo, Ontario N2L 3G1, Canada}}

\date{}

\maketitle
\begin{abstract}
Growing set of optimization and regression techniques, based upon sparse representations of signals, to build models from data sets has received widespread attention recently with the advent of compressed sensing.
{\red This paper deals with the parameterization of the Lorenz-96 (Lorenz, 1995) model with two time-scales that mimics mid-latitude atmospheric dynamics with microscopic convective processes.
Compressed sensing is used to build models (vector fields) to emulate the behavior of the fine-scale process, so that explicit simulations become an online benchmark for parameterization.
We apply compressed sensing, where the sparse recovery is achieved by constructing a sensing/dictionary matrix from ergodic samples generated by the Lorenz-96 atmospheric model, to parameterize the unresolved variables in terms of resolved variables. Stochastic parameterization is achieved by auto-regressive modelling of noise.
We utilize the ensemble Kalman filter for data assimilation, where observations~(direct measurements) are assimilated in the low-dimensional stochastic parameterized model to provide predictions. 
Finally, we compare the predictions of compressed sensing and Wilk’s polynomial regression  to demonstrate the potential effectiveness of the proposed methodology.}\\

\noindent
{\bf keywords:} compressed sensing, sparse regression, ensemble Kalman filter, auto-regression
\end{abstract}


\section{Introduction}
\label{S:Introduction}
Global weather and climate models are the best available tools for projecting likely climate changes.
However, current climate model projections are still considerably uncertain, due to uncertainties and errors in mathematical models, ``curse of dimensionality", model initialization, as well as inadequate observations and computational errors.
Multi-layer models are adequate to capture the dynamics of large-scale atmospheric phenomena.
Solving for small-scale processes with help of large eddy simulation~(LES) models, on a high-resolution grid embedded~(with grid spacing of a few tens of meters and can explicitly simulate the clouds and turbulence represented by parameterization schemes) within the low-resolution grid is very useful, but it's applicability remains limited to its high computational cost.
In this case, learning methods are used to develop parameterization schemes that utilize the explicitly simulated clouds and turbulence.
{\red In practice, low-resolution models resolve the large-scale processes while the effects of the small-scale processes are often parameterized in terms of the large-scale variables.}

In 1995, Ed Lorenz came up with a simple model that captured the essence of this problem with two stacked layers, which is referred to as the Lorenz-96~\cite{Lorenz1996} model. 
Even though this simple model is still a long way from the Primitive Equations that we need to eventually study, it is hoped that the parameterization schemes and filtering methods presented here can be ultimately adapted for the realistic models that need to be computationally solved for weather and climate predictions.
We use the Lorenz-96 model to demonstrate the sparse optimization methods and regression techniques to build models from data. This model, as discussed in Section~\ref{S:Lorenz96}, has nonlinearities resembling the advective nonlinearities of fluid dynamics and a multiscale coupling of slow and fast variables similar to what is seen in the two-layer atmospheric models.
{\red Resolving all the relevant length and time scales in simulations of the turbulent atmospheric circulations remains out of reach due to computational constraints.
Hence,
the small-scale processes must be represented indirectly, through parameterization schemes which estimate their net impact on the resolved atmospheric state.}

The objective of this paper is to develop a data-driven method to approximate the effects of small-scale processes with simplified functions called parameterizations which are characterized by compressed sensing (CS). 
Such parameterization schemes are inherently approximate, and the development of schemes which produce realistic simulations is a central challenge of model development.
Most parameterization schemes depend critically on various parameters whose values cannot be determined a priori but must instead be estimated through data.
Recent developments in optimization and sparse representations in high-dimensional spaces have provided powerful new techniques.

This paper deals with the parameterization of the Lorenz-96 model with two time-scale {\red simplified ordinary differential equations describing advection, damping and forcing}. 
We assume that the model that expresses the resolved atmospheric variables of the known physics or the resolved modes are well understood. 
Furthermore, Lorenz in introducing this multiscale model showed that the small-scale features of the atmosphere act like forcing.
{\red We present a CS and sparse recovery approach for constructing approximations to the Lorenz-96 model system that blends physics-informed candidate functions with functions that enforce sparse structure.}
The ultimate aim is to adapt the mathematical methods developed here for the realistic models that can be solved computationally to examine the turbulent atmospheric circulations. 
In the past few years, data-driven parameterization of the Lorenz-96 model using machine learning and deep learning techniques~\cite{Gagne2020} has shown promising results.   

{\red The aim of any parameterization is to develop enhanced models that capture the essence of dynamics and produce reliable predictions. Beside deterministic parameterization, stochastic parameterization was applied to the Lorenz-96 model with several kinds of noise and provided promising results \cite{Arnold13}. In our problem, we extract the useful elements for parameterization by using CS that relies on the restricted isometry property (RIP) \cite{Cande2006}. The RIP is satisfied due to the chaotic and ergodic nature of the Lorenz-96 model (numerically verified), showing that CS would be a great fit for stochastic parametrization. In this way, we are able to capture the key components that allow us to represent the unresolved variables in terms of resolved variables with only a few elements. After the parameterization of the deterministic part with compressed sensing, the stochastic part is achieved by auto-regressive modelling of noise. These parameterized models are then used for the predictions by using ensemble Kalman filter, where the  observations are assimilated through the slow time scale.}

The paper is organized as follows.
In Section~\ref{S:Lorenz96}, we will introduce the 40-dimensional Lorenz-96 model, that was originally suggested by \cite{Lorenz1996}, and describe the characterization of complex dynamics including statistical properties, Lyapunov exponents, marginal probability density functions (PDFs) and auto-correlation functions (ACFs).
Section~\ref{S:DDmethods} will present the details on data-driven reduction methods such as regression and compressed sensing with sparse optimization.
In Section~\ref{S:StochasticParameterization}, stochastic paramterization using compressed sensing approach will be illustrated on Lorenz-96 model.
Section~\ref{S:NonlinearFiltering} will describe a general ensemble filtering algorithm, where forecasting step is straightforward because the Markov kernel is known, while the update step is much more sophisticated. 
This nonlinear filtering method will be applied to parameterized models. 
Finally, conclusions and future work discussion will follow in Section~\ref{S:Conclusion}.
To ensure that our results are reproducible, we have provided our code in \cite{GithubCode}.

\section{The Lorenz-96 Atmospheric Model}
\label{S:Lorenz96}
\begin{figure}
\includegraphics[width=0.45\textwidth]{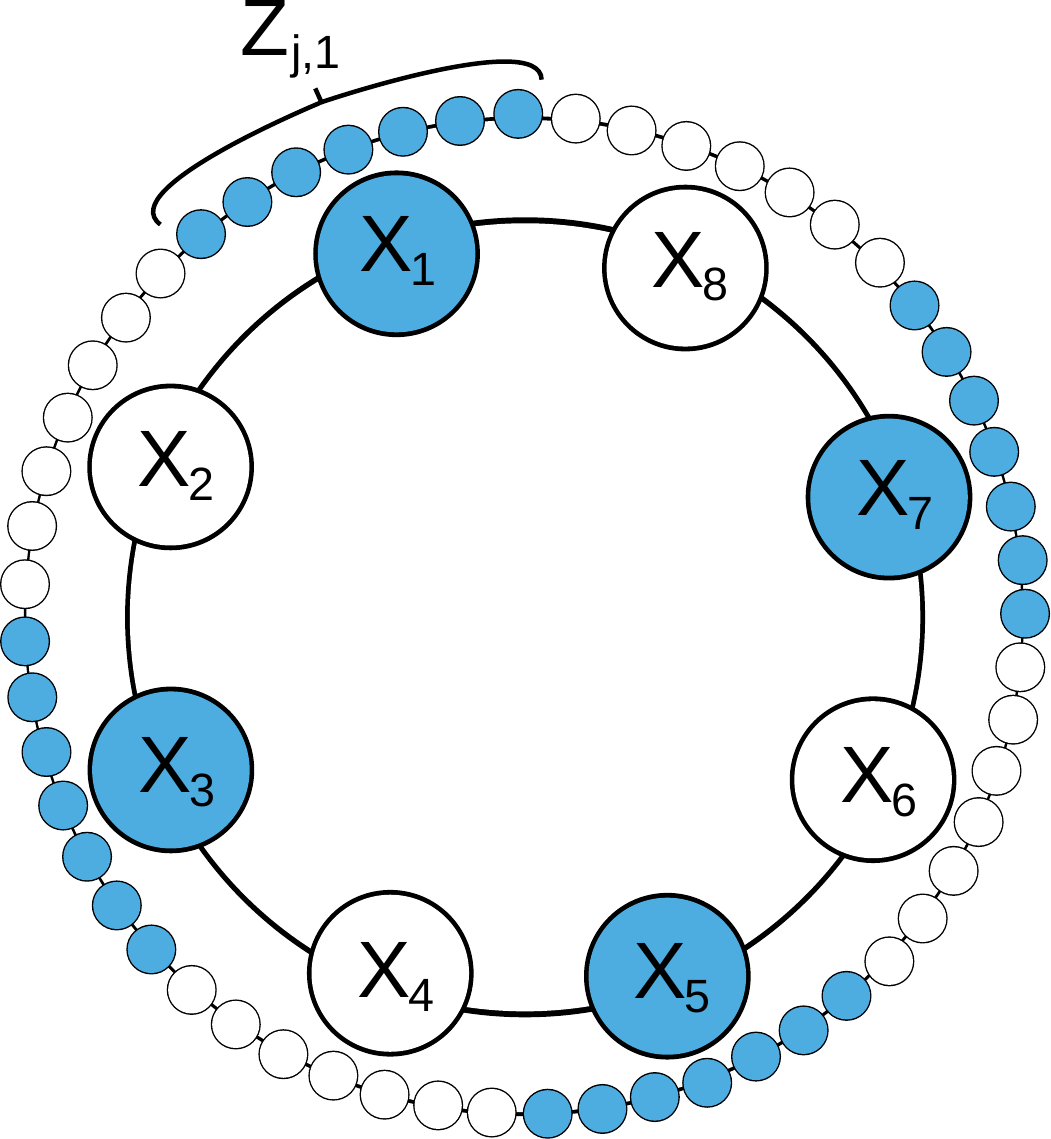}
\caption{A recreation of a similar image by Wilks (2005), that illustrates the two timescale Lorenz-96 model.}
   \label{fig:Lorenz}
\end{figure}

Let us now turn to a simple atmospheric model, which nonetheless exhibits many of the difficulties arising in realistic models, to gain insight into predictability and data assimilation. The Lorenz-96 model was originally introduced in~\cite{Lorenz1996} to mimic multiscale mid-latitude atmospheric dynamics for an unspecified scalar meteorological quantity at $K$ equidistant grid points along a latitude circle:
\begin{align}
\label{eq:2.0.1}\dot{X}_{k} &= - X_{k-1}(X_{k-2} - X_{k+1}) - X_{k} + F - \frac{hc}{b}\sum_{j=1}^{J} Z_{j,k} \\
\label{eq:2.0.2}\dot{Z}_{j,k} &= -cbZ_{j+1,k}(Z_{j+2,k} - Z_{j-1,k}) - cZ_{j,k} + \frac{hc}{b}X_{k}
\end{align}
Equation \eqref{eq:2.0.1} describes the dynamics of some atmospheric quantity $X$, and $X_k$ can represent the value of this variable at time $t$ in the $k^{\rm th}$ sector defined over a latitude circle in the mid-latitude region. 
({\it Note that we use subscripts $k$ and $j$ to conform with the typical spatial indexing notation used for the Lorenz-96 model. In sections that follow, subscripts $k$ and $j$ will be used as discrete time indices, not to be confused with the spatial indices of the Lorenz model}). 
The latitude circle is divided into $K$ sectors.\\
For convenience, the vector field associated with the slow components or the resolved variables is assumed a priori known with representable physics, which is denoted by 
\begin{equation}
\label{eq:prioriknown}
\G_{k}(X):=- X_{k-1}(X_{k-2} - X_{k+1}) - X_{k} + F 
\end{equation}
Figure \ref{fig:Lorenz} illustrates the behavior of a generic slow state $X_{k}$ (shown in the inner ring with $K  =8$ for illustration), the fast states in the $1$st sector, that is, $Z_{1,1}, \hdots, Z_{J,1}$ (shown in the outer ring with $J=7$ for illustration), and the fast scale forcing that enters~\eqref{eq:2.0.1}; for the $X_k$ component the fast scale forcing is 
\begin{equation}
\label{eq:forcing}
    U_k=-\frac{hc}{b}\sum_{j=1}^{J}Z_{j,k}
\end{equation}
Each $X_k$ is coupled to its neighbors $X_{k+1}$, $X_{k-1}$, and $X_{k-2}$ by \eqref{eq:2.0.1}. $X_k$ and $Z_{j,k}$ are periodic on the spatial domain, and \eqref{eq:2.0.1} and \eqref{eq:2.0.2} applies for all values of $k$ and $j$ by letting $X_{k+K} = X_{k}$, $Z_{j,k+K}=Z_{j,k}$, $Z_{j+J,k}=Z_{j,k+1}$ so, for example, in Figure \ref{fig:Lorenz} we have $X_{9} = X_{1}$, $Z_{1,9} =Z_{1,1}$, and $Z_{8,1} = Z_{1,2}$.
Even though this model is not a simplification of any atmospheric system, it is designed to satisfy three basic properties: {contains a quadratic discrete advection-like term that conserves the total energy (i.e. it does not contribute to $dE/dt$) just like many atmospheric models}; it has linear dissipation~(the --$X_{k}$ term) that decreases the energy defined as $E = \sum_{k=1}^{K}X_{k}$, and an external forcing term $F > 0$ that can increase or decrease the total energy. 

The original $K$ dimensional system was extended to study the influence of multiple spatio-temporal scales on the predictability of atmospheric flow by the dividing each segment $k$ into $J$ subsectors~($J=10$), and introducing a fast variable, $Z_{j,k}$ given by~\eqref{eq:2.0.2}, associated with each subsector. 
Thus, each $X_{k}$ represents a slowly-varying, large amplitude atmospheric quantity, with $J$ fast-varying, low amplitude, similarly coupled quantities, $Z_{j,k}$, associated with it.
In the context of climate modeling, the slow component is also known as the resolved climate modes while the rapidly-varying component is known as the unresolved non-climate modes.
The coupling terms between neighbors model advection between sectors and subsectors, while the coupling between each sector and its subsectors model damping within the model. 
The model is subjected to linear external forcing, $F$, on the slow timescale.

We consider $40$ dimensional slow variable components and $400$ dimensional fast variable components.
With $K=40$, one can think of $\mathbf{X}$ as representing an unspecified scalar meteorological quantity equally spaced along a mid-latitude belt of roughly 30,000 km, thereby giving the model the correct synodic scale to mimic Rossby waves. 
Our main research interest in this paper is to develop a dimensional reduction and parameterization method based on \textit{compressed sensing}~\cite{Cande2006,Donoho2006}. We compare this compressed sensing scheme based on sparse regression with~\textit{Wilk's polynomial regression}~\cite{Wilks2005}.

\subsection{Characterization of complex dynamics}
\label{S:ergodic}
According to~\cite{Devaney1980}, a chaotic dynamical system defined on a metric space $X$ by a continuous map $f$ is said to be chaotic on $X$ if it exhibits the following three essential features; (a) $f$ is transitive, (b) the periodic points of $f$ are dense in $X$, and (c) $f$ has sensitive dependence on initial conditions.
These systems are completely deterministic in that if the current state of the system $X_{0}$ is known, all future states at time $t$, denoted by $X_{t},t>0$, are known. However, even small numerical errors in the computation of $X_{0}$ can lead to very large differences in the values of $X_{t}$, so while these systems are deterministic, they are not predictable. Hence, chaotic systems exhibit a large departure in the long-term evolution for very tiny changes in initial conditions, and hence to study the limiting behavior of these systems, one cannot study the global evolution ``orbit by orbit". Instead, one studies the {\em statistical properties} of these systems, as we describe in this section.
\subsubsection{Statistical properties of chaotic systems}
\label{S:Statistical}
Although individual orbits of chaotic dynamical systems are by definition unpredictable, the average behavior of typical trajectories can often be given a precise statistical description. 
Anosov (1969), Pesin (2004) and Sinai (1977) showed that a wide class of ``chaotic" maps exhibit natural ergodic invariant measures, the maps being ``chaotic" in the sense that the orbits of two points, which may be arbitrarily close to each other with respect to a metric, may evolve very differently over time~(positive Lyapunov exponents as discussed below).
Loosely speaking, ergodicity is the ``law of large numbers" for dynamical systems. While the ergodic theorems prove convergence to the mean, natural questions arise about the deviation of these processes from the mean and about other statistical properties such as rates of mixing.
For a given invariant measure, and a class of observables, the correlation functions tell whether (and how fast) the system ``mixes", i.e. ``forgets" its initial conditions. 

Many statistical properties, such as the Central Limit Theorem and the Law of Iterated Logarithm, displayed by independent identically distributed stochastic processes, have been shown to hold for stochastic processes generated by chaotic dynamical systems.
One of the better known techniques for establishing such results was proved by Young (1998), where it was shown that for such~(chaotic) dynamical systems with certain ``towers" property~(a powerful tool for coding the dynamics of a system in such a way so
as to enable the computation of its statistical properties), behave as Markov chains on spaces with infinitely many states.
Should such dynamical systems mix sufficiently quickly, the process satisfies the {\em central limit theorem} and the completely deterministic process behaves, at least in the first two moments, as an independent stochastic process. Furthermore, the degree of independence is characterized by the rate at which the process mixes.
\subsubsection{Lyapunov Exponents (LEs)}
\label{S:LEs}
Among the features listed above, sensitivity to initial conditions is widely known to be associated with chaotic dynamical systems, and a quantitative measure of this sensitivity is given by the ${\it Lyapunov\; exponents}$ (LEs) which are the most commonly used quantities to characterize chaotic phenomena.
The LEs are asymptotic measures characterizing the average rate of divergence (or convergence) of small perturbations to the solutions of a dynamical system. The value of the maximum LE is an indicator of the chaotic or regular nature of orbits.

Considering the phase space, which is assumed to be $\mathbb{R}^{d}$ or a Riemannian manifold, and denoting by $M$ throughout, the Lebesgue or the Riemannian measure on $M$ is denoted by $m$, and the dynamics are generated by iterating a self-map of $M$, written $f : M \circlearrowleft$. Given a differentiable map $f: M \circlearrowleft$, a point $x \in M$ and a tangent vector $v$ at $x$, we define
\begin{equation}
\lambda (x, v) = \lim_{n \rightarrow \infty} \frac{1}{n} \log |Df_{x}^{n} (v)| \hspace{0.5cm} \text{if this limit exists}
\end{equation} 
{\red where $Df_{x}^{n} (v)$ is the derivative of  the $n$th iteration of $f$, (i.e., $f\circ \cdots \circ f$), with respect to variable $x$ at the vector $v$.}
And if the limit does not exist, then the limit is replaced by $\liminf$ or $\limsup$ and we write        $\underline{\lambda}(x,v)$ and $\overline{\lambda}(x,v)$, respectively. Thus, $\lambda(x,v) > 0$ implies the exponential growth of $|Df_{x}^{n} (v)|$  and can be interpreted to mean exponential divergence of nearby orbits.
The computational procedure for Lyapunov exponents using the standard method can be found in~\cite{Wolf1985,Skokos2010}.

When the chaotic system is dissipative, solutions will eventually reside on an attractor.
The dimension of this attractor, i.e., the number of ``active" degrees of freedom, is also known as the Kaplan-Yorke dimension 
$D_{KY}$~\cite{Eckmann1985} and is given as
\begin{equation}
    D_{KY} = r + \frac{1}{|\tilde{\lambda}_{r+1}|}\sum_{i=1}^{r}\tilde{\lambda}_{i}
\end{equation}
where $\tilde{\lambda}_1 \ge \tilde{\lambda}_2 \ge \cdots \ge \tilde{\lambda}_d$ are the Lyapunov exponents counted with multiplicity ($d = \mbox{dim}(M)$) and $r$ is the largest integer such that $\sum_{i=1}^r \tilde{\lambda}_{i} > 0$.

The Kolmogorov-Sinai entropy $H_{KS}$ is a measure of the diversity of the trajectories generated by the dynamical system. 
This entropy measure is determined through the Pesin formula, which provides an upper bound to $H_{KS}$~\cite{Eckmann1985} as
\begin{equation}
    H_{KS} = \sum_{\tilde{\lambda}_i > 0} \tilde{\lambda}_{i}.
\end{equation}
Complexity of the dynamics of the Lorenz-96 model varies considerably with different values of the constant forcing term $F$. 
It is shown in~\cite{Majda2005}, the model is in weakly chaotic regimes for $F=5,\,6$, strongly chaotic regime for the values of $F$ from $8$ to $10$, and turbulent regimes for $F=12,\,16,\,24$.
 We will present our results for a specific value of forcing $F=10$.

\begin{center}
\begin{figure}
\includegraphics[width=0.50\textwidth]{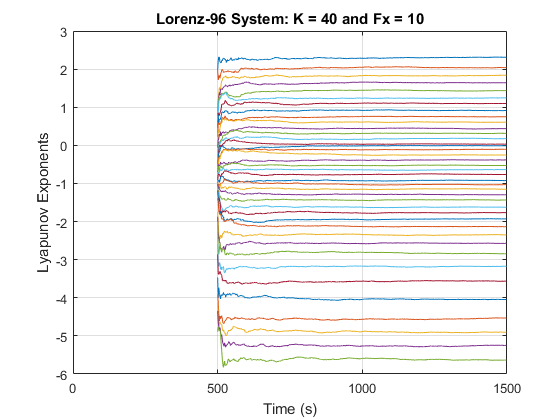}
    \caption{Time evolution of Lyapunov exponents for the Lorenz-96 system}
    \label{fig:LE_traj}
\end{figure}
\end{center}

In our numerical simulations, we begin from random initial conditions and use a fourth-order Runge-Kutta time integration with a time step of 0.01. 
Typically, we integrate forward for 500 time units before starting the calculation of Lyapunov exponents to ensure that all transients have decayed.
At this point we begin integrating the tangent space equations using the Gram-Schmidt procedure update every $\tau=0.2$ time units.
With $K = 40$, $F = 10$, $h=0$~(just the slow variables), $c=10$ and $b=10$, the summary of Lyapunov Exponents obtained from the above simulation of the Lorenz-96 system is given Table 1. 
The time evolution of the Lyapunov exponents for the Lorenz-96 model is shown in Figure \ref{fig:LE_traj}. As noted before,the calculation of the Lyapunov exponents begins at 500 time units and Figure \ref{fig:LE_traj} depicts the convergence of the Lyapunov exponents around 1000 time units.
\begin{table}
\centering
\begin{tabular}{ |l|l| } 
 \hline
 Largest Lyapunov exponent $\lambda_1$ & 2.3098 \\ 
 Error doubling time &  0.3 time units \\ 
 Number of strictly positive Lyapunov exponents &  14 \\ 
 Number of neutral, $\lambda \in [- 1E-2, 1E-2]$, exponents & 1  \\
 Number of strictly negative Lyapunov exponents &  26 \\
 Kaplan-Yorke dimension   &  29.4694  \\
 Kolmogorov-entropy   &  14.8409  \\
 \hline
\end{tabular}
\caption{Summary of Lyapunov Exponents}
\label{table:LyapunovExponents}
\end{table}
\subsubsection{Marginal Probability Density Functions (PDFs) and Auto-correlation Functions (ACFs)}
Typically, the PDFs are computed using bin-counting also known as histogram methods. For estimating the PDF of a continuous variable $x$, standard histogram methods simply partition $x$ into distinct bins of width $\Delta_{i}$, and then the number of observations $n_i$ of $x$ falling into bin $i$ is converted into a normalized probability density for that bin.

\begin{figure}
\includegraphics[width=0.50\textwidth]{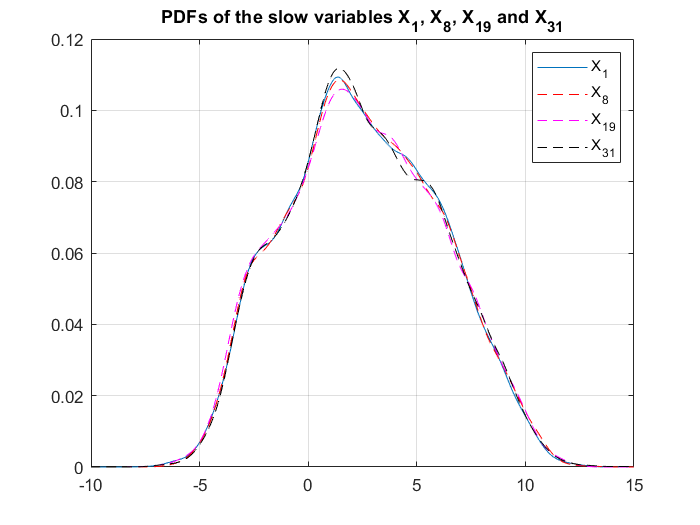}
    \caption{PDFs of the slow variables $X_1$, $X_8$, $X_{19}$ and $X_{31}$}
    \label{fig:PDF-X1}
\end{figure}

More smoother estimates for PDFs can be obtained by using the Gaussian product kernel function which results in the estimation of the PDF known as the {\it kernel density estimation} (KDE)~\cite{Bishop2006}:
\begin{equation}
    \label{eq:KernelDensityEstimation}
    p_{X}(x) = \frac{1}{N}\sum_{j=1}^{N}\frac{1}{(2\pi\lambda^2)^{1/2}}\exp\{-\frac{\|x_{j}-x\|^2}{2\lambda^2}\}
\end{equation}
where $\lambda$ represents the standard deviation of the Gaussian components and $N$ is the total number of data points.
Figure \ref{fig:PDF-X1} depicts the marginal PDFs of some (randomly selected) slow variables $X_1$, $X_8$, $X_{19}$ and $X_{31}$ simulated as described above in Section 2.
The convergence of the marginal PDFs shown in Figure \ref{fig:PDF-X1} demonstrates the existing symmetry among the slow variables.

\begin{figure}
\includegraphics[width=0.50\textwidth]{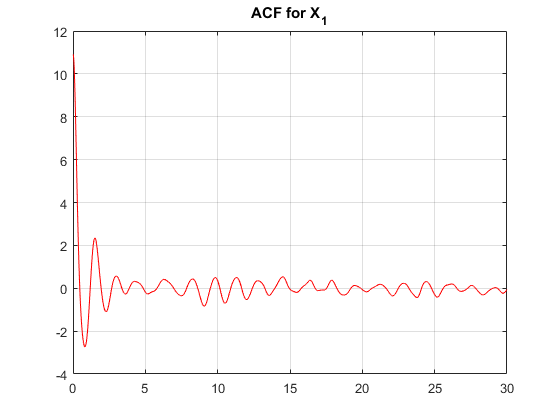}
\caption{ACF of the slow variable $X_1$}
\label{fig:ACF-X1}
\end{figure}

The strongly mixing character of the slow dynamics can be inferred from the correlation function, which means that the correlation function decays quickly to zero or near zero values. 
The discrete version of the auto-correlation functions (ACFs) for the slow variables is given by
\begin{align}
& C_{k,k}(m\Delta t) = \nonumber \\ &\frac{1}{M-m}\sum_{h=1}^{M-m}(X_{k}(h\Delta t)-\bar{X}) (X_{k}((h+m)\Delta t) - \bar{X})
\end{align}
where
\begin{equation}
\bar{X} = \frac{1}{M}\sum_{m=1}^{M}X_{k}(m\Delta t)
\end{equation}
and these formulae can be given similarly for the fast variables. Figure \ref{fig:ACF-X1} shows the ACF for the slow variable $X_1$ which decays very quickly and
then oscillates around zero illustrating the strongly mixing character of the dynamics.

\section{Data-driven Reduction Methods}
\label{S:DDmethods}
The objective of this section is to compare model-based and data-driven tools for dimensional reduction and reducing uncertainties in prediction.
It is impossible to effectively ``learn" from high dimensional data unless there is some kind of implicit or explicit low-dimensional structure -- this can be made mathematically precise in different ways. The methods presented below provide the low-dimensional structure either because of the presence of multi-scales in the model or due to the measurement matrix chosen to satisfy with high probability the so called ``Restricted Isometry Property (RIP)"  condition introduced in compressed sensing~\cite{Cande2006}. 

Compressed sensing theory addresses the problem of recovering high-dimensional but sparse vectors from a limited number of measurements - 
 the number of observed measurements $m$ is significantly smaller than the dimension $n$ of the original vector.
{\red We present a compressed sensing approach for sparse regression problem and using the ideas from concentration of measure, extend it to a case where the samples are generated from a chaotic time series, which are sufficiently mixing, thus ergodic, as shown in subsection~\ref{S:ergodic}. Here, the sparse recovery is achieved by constructing a matrix $\Theta $ in section~\ref{S:CS} from ergodic samples generated by the Lorenz-96, such that it is suitable for reconstructing sparse parameter vectors via $l_1$-minimization.}

\subsection{Regression}
\label{S:Regression}
This parameterization was proposed by Wilks~\cite{Wilks2005}, and it uses a polynomial regression to parametrize the effects~\eqref{eq:forcing} of unresolved variables:
\begin{equation}
    U_k(t) = P(X_k(t)) + e_k(t)
\end{equation}
where $P$ is a 4th degree polynomial in $X_k$ and is given as:
\begin{equation}
    P(X_k) = a_0 + a_1X_k + a_2X_k^2 + a_3X_k^3 + a_4X_k^4
\end{equation}
And $e_k$ represents noise, which can be written as a first-order autoregressive model:
\begin{equation}
    e_{k}(t_i)=\phi e_{k}(t_{i-1}) + \sigma z_{k,i}
    \label{eq:ARMA}
\end{equation}
where $z_{k,i}\stackrel{i.i.d}{\sim}N(0,1)$ for all $i,k$. The method of calculating $\phi$ and $\sigma$ is explained in Appendix A.\\

\subsection{Compressed Sensing and Sparse Recovery}
\label{S:CS}
Over the past two decades, researchers have focused on sparsity as one type of low-dimensional structure. 
Given the recent advances in both compressed sensing~\cite{Cande2006,Donoho2006,Cande2008} and sparse regression~\cite{Tibshirani1996} it has become computationally feasible to extract system dynamics from large, multimodal datasets. 
The term sparse in signal processing context refers to the case where signals~(or any type of data, in general) have merely few non-zero components with respect to the total number of components.
Sparsity plays a key role in optimization and data sciences.
In the context of our work, these techniques rely heavily on the fact that many dynamical systems can be represented by governing equations that are sparse in the space of all possible functions of a given algebraic structure.

{\red Compressed sensing~(CS) is a technique for sampling and reconstructing signals as sparse signals, that is, the reconstructed signals have few non-zero components.}
More precisely, $k$-sparse signals are those that can be represented by $k<<n$ significant coefficients over an $n$-dimensional basis. 
The central goal of CS  is to capture attributes of a signal using very few measurements, which distinguishes CS from other dimensionality reduction techniques - a methodology for the recovery of sparse vectors from a small number of linear measurements. Hence, this allows for polynomial-time reconstruction of the sparse signal~\cite{Donoho2006}.

{\red In Donoho (2006) and Cand\`{e}s and Tao (2006), the original signal $X$ is projected onto a lower-dimensional subspace via a random projection scheme, called the sensing matrix.}
More precisely, this broader objective is exemplified by the important special case in which one is interested in finding a vector {\red $S\in\R^n$} using  the 
(noisy) observation or the measurement data 
\begin{equation}
\red Y=\Theta S +\eta
\end{equation}
{\red where $ \Theta~(=:\Theta(X))~\in\C^{m\times n}$ with $k<m<n$ is the dictionary or sensing matrix and $\eta$ is the measurement noise. In CS, we estimate a good sensing matrix $\Theta$ and a vector $S$ that reconstructs the signal $Y$ from $X$.}

In general, this problem cannot be solved uniquely. However, if {\red $S$} is $k$-sparse i.e., if it has up to $k$ nonzero entries, the theory of CS shows that it is possible to reconstruct $X$ uniquely, from $m$ linear measurements even when $m<<n$, by exploiting the sparsity of {\red $S$}. This can be achieved
by finding the sparsest signal consistent with the vector of measurements~\cite{Donoho2006}, i.e.
\begin{equation}
\label{P:l0}\red
\arg\;\min_{S\in \R^{n}}
\|S\|_{0} \quad \text{subject to}\quad \|Y-\Theta\,S\|_{_{2}}\leq \varepsilon,
\end{equation}
where {\red $\|S\|_{0}$} denotes the $l_{0}$ norm for {\red $S$} (i.e., the number of non-zero entries of $S$), while $\varepsilon$ denotes a parameter that depends on the level of measurement noise $\eta$.
It can be shown that the $l_{0}$  minimization method can exactly reconstruct the original signal in the absence of noise using a properly chosen sensing matrix $\Theta$ 
whenever $m>2k$. However, $l_{0}$ minimization problem~\eqref{P:l0} is non-convex and combinatorial, which is NP-hard.

Hence, instead of problem~\eqref{P:l0}  we consider its $l_{1}$ convex relaxation which may be stated as~\cite{Chen1998}
\begin{equation}
\label{P:l1}
\red
\arg\;\min_{S\in \R^{n}}
\|S\|_{1} \quad \text{subject to}\quad \|Y-\Theta\,S\|_{_{2}}\leq \varepsilon,
\end{equation}
where $\|\cdot\|_{1}$ represents the $l_1$ norm, which is a convex function and problem~~\eqref{P:l1} is a convex optimization problem which can accurately approximate the signal $X$~(solution to problem~\eqref{P:l0}) in polynomial time if the sensing matrix $\Theta$ is chosen to satisfy the necessary RIP condition with high probability~\cite{Cande2005,Cande2006}.
Loosely speaking, if $\Theta$ satisfies the RIP, then the measurement matrix approximately preserves the Euclidean length of every signal. Equivalently, all subsets of $k$ columns taken from $\Theta$  are nearly orthogonal.
One should note that the $l_{1}$ minimization in \eqref{P:l1} is closely related to
the Lasso problem~\cite{Tibshirani1996}.
\begin{equation}
\label{P:l2}
\red
  \arg\;\min_{S\in \R^{n}}
  \frac{1}{2}\|Y-\Theta\,S\|_{_{2}}^2 +\lambda\,\|S\|_{1}
\end{equation}
where $\lambda \geq 0$ is a regularization parameter. If $\varepsilon$ and $\lambda$ in~\eqref{P:l1} and~\eqref{P:l2} satisfy some special conditions, the two problems are equivalent; however, characterizing the relationships between $\varepsilon$ and $\lambda$ is difficult, except for the special case of orthogonal sensing matrices $\Theta$.

Even though there are several objectives in compressed sensing, we are mainly interested in {\red estimating a good sensing matrix and} efficiently reconstructing the $k$-sparse vector {\red $S$} from the measurements. Hence, the RIP plays an important role in our calculations.
Unfortunately, direct verification of RIP for a given matrix is not feasible. 
Nevertheless, it has been shown that there are some matrices, particularly random matrices with independent and identically distributed (i.i.d) entries, satisfy RIP with high probability.

We now enforce this sparse structure within learning. 
One of the related problems in the context of sparsity is ``dictionary learning".
We refer to the columns of the ${m\times n}$ matrix  $\Theta$, as a dictionary, and they are assumed to be a set of vectors capable of providing a highly succinct representation for the signal. 
One has freedom to determine the particular dictionary elements; however, the choice of monomial terms provides a model for the interactions between each of the components of the time-series and is used for model identification of dynamical systems in the general setting.
In  dictionary  learning, which we shall use to determine $\Theta$, the aim is to find a dictionary that can sparsely represent a given set of training data for further analysis --``sparsifying matrix" from a given set of training data. 
Regularized optimization with polynomial dictionaries is used to approximate the generating function of some unknown dynamic process.
When the dictionary is large enough so that the ``true" function is in the span of the candidate space, the solutions produced by sparse optimization are guaranteed to be accurate.
This means that any of the training data should be presentable using linear combinations of few columns from the dictionary. 
{\red In Lorenz-96, the  sensing matrix $\Theta$ is formed as the collection of matrix of all monomials (stored column-wise) of resolved modes and due to the ergodic nature of the simulated time series the entries satisfy RIP with high probability.}

In this paper, the effects of unresolved processes~\eqref{eq:2.0.2} on the resolved modes~\eqref{eq:2.0.1} in the Lorenz-96 model are denoted by $u_k=- \frac{hc}{b}\sum_{j=1}^{J} Z_{j,k}$.
The CS method takes snapshot data ${\bf x}(t)\in\R^n$ and attempts to discover the structure of a potentially highly-nonlinear mapping of the form
\begin{equation}
{\bf u}(t) \approx {\bf f}({\bf x}(t))
\end{equation}
where the $n$-dimensional vector ${\bf x}(t)=[x_{1}(t) \ x_{2}(t) \ \hdots \ x_{n}(t)]^{T}$ represents the state of the system at time $t$.
{\red The error ${\bf e}(t) = {\bf u}(t) - {\bf f}({\bf x}(t))$ at time $t$ will be used later for auto-regression.}
The basis of this construction lies on the fact that in the multi scale Lorenz-96 model, the unresolved dynamics depends on the value of the resolved variables $\bf x$.

In particular, this method searches for a sparse approximation of the function mapping ${\bf f}({\bf x})=[f_{1}({\bf x}) \ f_{2}({\bf x}) \ \hdots \ f_{n}({\bf x})]^{T}\in\R^n$
using a dictionary of basis functions yielding
\begin{equation}
f_{k}({\bf x}) \approx \sum_{j = 1}^{p} \theta_j({\bf x})s_{jk} = {\bf \Theta}({\bf x}){\bf s}_k
\end{equation}
where ${\bf \Theta}({\bf x})=[\theta_{1}({\bf x}), \theta_{2}({\bf x}),...,\theta_{p}({\bf x})]\in\R^{1\times p}$ form the dictionary of basis functions, and ${\bf s}_k=[s_{1k} \ s_{2k} \ \hdots \ s_{pk}]^{T}\in\R^p$ is the vector of coefficients.
Sparsity is expressed ``not" in terms of an orthonormal basis but in terms of an over-complete dictionary.
The majority of the coefficients $s_{jk}$ are zero while the remaining nonzero entries identify the active terms contributing to the sparse representation of ${\bf f}({\bf x})$.

To learn the function mapping $\bf{f}$ from data, time-histories of vectors ${\bf x}(t)$ and ${\bf u}(t)$ are collected by performing simulations for the Lorenz-96 system.
The input data matrix ${\bf X}\in\R^{m\times n}$ is made of $m$ snapshots of the state vector sampled at several times $t_{1}, t_{2},...,t_{m}$ where subscripts $1,2,\hdots ,n$ indicates the elements of the state vector. The output data matrix ${\bf U}\in\R^{m\times n}$ follows the same arrangement, and these two data matrices are given as follows:
\begin{eqnarray*}
\bf{X} = \left[
\begin{array}{c}
{\bf x}^{T}(t_1) \\
{\bf x}^{T}(t_2) \\
\vdots \\
{\bf x}^{T}(t_m)
\end{array}
\right] = \left[
\begin{array}{cccc}
x_1(t_1) & x_2(t_1) & \hdots & x_n(t_1) \\
x_1(t_2) & x_2(t_2) & \hdots & x_n(t_2) \\
\vdots & \vdots & \ddots & \vdots \\
x_1(t_m) & x_2(t_m) & \hdots & x_n(t_m)
\end{array}
\right],
\end{eqnarray*}
\begin{eqnarray*}
{\bf U} = \left[
\begin{array}{c}
{\bf u}^{T}(t_1) \\
{\bf u}^{T}(t_2) \\
\vdots \\
{\bf u}^{T}(t_m)
\end{array}
\right] = \left[
\begin{array}{cccc}
u_1(t_1) & u_2(t_1) & \hdots & u_n(t_1) \\
u_1(t_2) & u_2(t_2) & \hdots & u_n(t_2) \\
\vdots & \vdots & \ddots & \vdots \\
u_1(t_m) & u_2(t_m) & \hdots & u_n(t_m)
\end{array}
\right].
\end{eqnarray*}

\subsubsection{Dictionary of Basis Functions}
\label{S:Dictionary}
First, we construct an augmented dictionary 
matrix ${\bf \Theta}({\bf X})\in\R^{m\times p}$ where each column represents a basis function that can be a potential candidate for the terms in ${\bf f}$ to be discovered from data. $p$ is the maximal number of $n$-multivariate monomials of degree at most $d$. As stated in~\cite{Brunton2016}, there is huge flexibility in choosing the basis functions to populate the dictionary 
matrix ${\bf \Theta}({\bf X})$.  For example, the dictionary 
matrix ${\bf \Theta}({\bf X})$ may consist of constant, polynomial and trigonometric terms as shown in~\eqref{eq:3.2.1.1}.
%
\begin{eqnarray}
{\bf \Theta}({\bf X}) = \left[
\begin{array}{ccccccccc}
| & | & | & \hdots & | \\
\bf{1} & \bf{X} & \bf{X}^{P_2}  & \hdots & \bf{X}^{P_d} \\
| & | & | & \hdots & | 
\end{array}
\right]_{m\times p}
\label{eq:3.2.1.1}
\end{eqnarray}
%
The dictionary matrix ${\bf \Theta}({\bf X})$ is constructed by stacking together, column wise, candidate nonlinear functions of ${\bf X}$. 
Here, higher order polynomials are denoted as ${\bf X}^{P_d}$ where $d$ is the order of the polynomial considered. 
For example, element 1 is a column-vector of ones, element ${\bf X}$ is as defined above, element ${\bf X}^{P_2}$ is the matrix containing the set of all quadratic polynomial functions of the state vector ${\bf x}$, and is constructed as shown in~\eqref{eq:3.2.1.2}.
\begin{figure}
\begin{eqnarray}
\bf{X}^{P_2} = \left[
\begin{array}{cccccc}
x_1^{2}(t_1) & x_1(t_1)x_2(t_1) &\hdots & x_2^{2}(t_1) & \hdots & x_n^{2}(t_1) \\
x_1^{2}(t_2) & x_1(t_2)x_2(t_2) &\hdots & x_2^{2}(t_2) & \hdots & x_n^{2}(t_2) \\
\vdots & \vdots & \ddots & \vdots & \ddots & \vdots \\
x_1^{2}(t_m) & x_1(t_m)x_2(t_m) & \hdots & x_2^{2}(t_m) & \hdots & x_n^{2}(t_m)
\end{array}
\right]
\label{eq:3.2.1.2}
\end{eqnarray}
\end{figure}

\subsubsection{Sparse Optimization}
\label{S:SparseOptimization}
Sparse regression is performed to approximately solve 
\begin{equation}
{\bf U} \approx {\bf \Theta}({\bf X}){\bf S}
\label{eq:lsqprob}
\end{equation} 
to determine the vector coefficients ${\bf S} = [{\bf s}_1 \ {\bf s}_2 \ \hdots \ {\bf s}_n]\in\R^{p\times n}$ that determines 
the active terms in ${\bf f}$. 
Since only a few of the candidate functions
are expected to have an effect on the function mapping ${\bf f}$, all coefficient vectors ${\bf s}_i$  are expected to be sparse.
In this work, we use the Lasso algorithm~\cite{Tibshirani1996} for solving the above sparse regression problem. 

The goal of the sparse optimization to solve the regression problem~\eqref{eq:lsqprob} can be defined using the Lasso form~\eqref{P:l2} for $k=1,...,n$
\begin{equation}
{\bf s}_{k} =\underset{{\bf s}_{k}^{'}}{argmin} \; ||{\bf u}_{k} - {\bf \Theta}({\bf X}){\bf s}_{k}^{'}||_{2}^2 + \lambda ||{\bf s}_{k}^{'}||_{1}
\end{equation}
which uses $l_1$ norm constraint on the coefficient vector ${\bf s}_k$.
The term $\lambda ||{\bf s}_{k}^{'}||_{1}$ is the regularization term, which
penalizes coefficients different from 0 in a linear fashion. It is
the actual promoter of sparsity in the minimization problem.

In the Lasso method, given a collection of $m$ time snapshot pairs $({\bf x}(t_i), {\bf u}(t_i))_{i=1}^{m}$, we estimate the coefficients in \eqref{eq:lsqprob} by solving the following optimization problem for $k = 1,\hdots,n$
\begin{equation}
{\red \underset{{\bf s}_{k}}{minimize} \; \sum_{i=1}^{m}\left(u_{k}(t_i) - \sum_{j=1}^{p}\theta_{j}({\bf x}(t_i))s_{jk}\right)^2 + \lambda \sum_{j=1}^{p}|s_{jk}|}
\end{equation}
{\red for some Lagrange multiplier $\lambda \geq 0$. In this paper, we use $\lambda=10^{-3}$.}

\section{Stochastic Parameterization: Application to Lorenz-96}
\label{S:StochasticParameterization}
Having introduced the data-driven methods for model reduction in the presence of time-scale separation, let us now turn to a simple atmospheric model, which nonetheless exhibits many of the difficulties arising in realistic models, to gain insight into predictability and data assimilation.
The dynamics of the Lorenz-96 equation are shown in~\eqref{eq:2.0.1} and \eqref{eq:2.0.2}. Let $X_k$ and $Z_{j,k}$ be the numerical solutions of the Lorenz-96 equations, generated using the 4th order Runge-Kutta method, for $k=1,...,K$, $j=1,...,J$. We define $U_k$ as in Equation \eqref{eq:forcing}:
\begin{equation*}
    U_k=-\frac{hc}{b}\sum_{j=1}^{J}Z_{j,k}
\end{equation*}
The resolved modes in~\eqref{eq:2.0.1} is now modified as:
\begin{equation}
    \label{eq:4.2.1}\dot{X}_{k} = - X_{k-1}(X_{k-2} - X_{k+1}) - X_{k} + F + U_k= \G_{k}(X) + U_k
\end{equation}
Our goal is to use the methods discussed in Section~\ref{S:DDmethods} to form parameterizations, in other words, express $U_k$ as a function of $X_1,...,X_{K}$ and compare each of them. Based on the ergodic properties of the Lorenz-96 system shown in Section~\ref{S:ergodic}, the Lyapunov exponents stabilize at $t=500$. The training set is set to $t_{train}=[500, 500+h, 500+2h, ..., 1000]$ where $h=0.01$. The test set is set to $t_{test}=[1000, 1000+h, 1000+2h, ..., 2000]$.

Using polynomial regression to solve stochastic parameterizations in the Lorenz-96 system was proposed by Wilks (2005). The goal of this method is to find a degree-4 polynomial $P$ such that $U_k$ can be approximated by $P(X_k)$ for $k=1,...,K$. This requires solving a least squares problem:
\begin{equation}
    \sum_{k=1}^{K}\sum_{i\in t_{train}}(U_k(t_i)-P(X_k(t_i)))^2
\end{equation}
The polynomial obtained from this method is:
\begin{align}
 \label{WilksPolynomial}
    P(X_k) &= -0.0003142 X_k^4 + 0.004882 X_k^3  \nonumber \\& + 0.005197 X_k^2 - 0.47362 X_k - 0.03779
\end{align}

The goal of compressed sensing in this context is to find a sparse matrix $S$ that solves for $Y = \Theta S$. Let $Y$ be the matrix with columns $U_1,...,U_K$, each divided by its $L_2$ norm. Every column of $\Theta$ is a basis function of $X_1,...,X_k$ divided by its $L_2$ norm. The goal of compressed sensing is to find unique equations for every $U_k$ for $k=1,...,K$. Since regression gave us a polynomial of degree 4, our goal is to express $U_k$ as a function of $X$ values of order 4. However, with $K=40$, this means the $\Theta$ matrix will have $135750$ columns, and thus, compressed sensing will be computationally expensive. Each of the parameterizations we get are shown in the "CompressedSensingEquations" folder in  \cite{GithubCode}.

A computationally less expensive approach will be to find a solution of order 2. 
{\red The 40-dimensional CS model is used for comparison of distributions, but we present below the single equation with averaged coefficients that will be used later for filtering}:
\begin{equation}
    \label{eq:4.2.2}{\red f_k(\textbf{X})}=-0.376X_k+0.0166X_k^2
\end{equation}
This parameterization shows that every $U_k$ is expressed as a function of $X_k$ and $X_k^2$, thus showing that every basis function $X_iX_j$ where $i\neq j$ is not relevant, and the $\Theta$ matrix can be reduced to a matrix with $X_k$ and $X_k^2$ as the columns.

We then explore solutions of order 3. Given $K=40$, $X^{P_3}$ has $11480$ columns, thus finding a solution of order 3 is computationally expensive. However, creating multiple parameterizations, each with a subset of columns of $X^{P_3}$ is computationally less expensive. Thus, in this section, we create 100 parameterizations, each of which the $\Theta$ matrix contains $X_k^3$ and 100 random columns from $X^{P_3}$. If we average the coefficients of the model, we get:
\begin{equation}
    \label{eq:4.2.3}{\red f_k(\textbf{X})}=-0.0338X_k+0.00452X_k^2+0.00142X_k^3
\end{equation}

This shows that every $U_k$ is expressed as a function of $X_k$, $X_k^2$ and $X_k^3$, thus showing that all remaining columns in $X^{P_3}$ are not relevant. Based on the results of the previous models, we expect $X_1^4,...,X_K^4$ to be the only relevant columns from $X^{P_4}$. {\red This information captures the underlying structure in the Lorenz-96 equation since, in equation (\ref{eq:2.0.2}), $\dot{Z}_{j,k}$ only depends on $X_k$.}

Thus, our $\Theta$ matrix contains the basis functions $1,X_k,X_k^2,X_k^3,X_k^4$ for $k=1,...,K$. We will call this parameterization the raw degree-4 compressed sensing model. Since this algorithm will be used for testing purposes, the compressed sensing algorithm also calculates biases for each $U_k$. If we average the coefficients, we get:
\begin{align}
    \label{eq:4.2.4}{\red f_k(\textbf{X})}&=-0.416X_k+0.00277X_k^2  \nonumber \\&+0.00509X_k^3-0.000306X_k^4-0.0243
\end{align}
The biases are calculated by the compressed sensing algorithm as $[2.173746, 1.393404, ..., -2.973427]$ and their mean is $-0.0243$. The problem with this parameterization is that there is too much variation in the biases. 

We investigate other ways to modify the biases. In the following model, the biases are strictly set to zero while keeping other coefficients constant. If we average the coefficients, we get:
\begin{align}
    \label{eq:4.2.5}{\red f_k(\textbf{X})}&=-0.474X_k+0.00277X_k^2  \nonumber \\
    &+0.00509X_k^3-0.000306X_k^4
\end{align}
In the next model, the biases are strictly set to the average of the biases while keeping other coefficients constant. If we average the coefficients, we get the same equation as (\ref{eq:4.2.4}). Both of these models perform equally well in predicting the trajectory of the system and significantly better than the raw 4-dimensional model. The model where biases are strictly set to the average of the biases performs better in the Kullback-Leibler divergence measure, but does not outperform linear regression.

A small variability was introduced in the biases with the goal of outperforming linear regression in terms of the Kullback-Leibler divergence measure. The biases were selected by sampling from $N(\mu,\sigma^2)$, where $\mu$ is the average of the biases ($\mu=-0.0243$). {\red We used $\sigma=0.07$ as it gave the lowest prediction error on the training set.} This leads to a new model, which outperforms all other models. The model is shown in Appendix B. 
If we average the coefficients, we get:
\begin{align}
    {\red f_k(\textbf{X})} &= -0.474X_k+0.00277X_k^2  \nonumber \\
    &+0.00509X_k^3-0.000306X_k^4-0.0241
    \label{eq:FinalCSModel}
\end{align}
We will use (\ref{eq:FinalCSModel}) as our compressed sensing model in the later sections.

The equations of the resolved variables associated with the vector field denoted by $\G_{k}(X)$ are equivariant with respect to a cyclic permutation of the variables. 
Therefore, the system with dimension $K$ is completely determined by the equation for the $k^{\rm th}$ variable. 
Even though the unresolved dynamics depend on the value of the resolved variables $\bf X$, while computing the parameterization of $U_{k}$ by CS, there were no restrictions placed on the form of the polynomial on the sparse regression. 
The consequences of ${\bf \Theta}({\bf X})$ satisfying RIP~(hence, sparsity)  are quite profound as shown in results provided in Appendix B. Compared to regression, compressed sensing also has the advantage that each unresolved variable $U_k$ is given a unique equation that best represents its dynamics as a function of resolved variables.

\subsection{Comparison of Trajectories}
\label{S:Trajectories}
Section~\ref{S:LEs} shows that the largest Lyapynov exponent is 2.3. This shows that the error-doubling time is $\frac{\ln(2)}{2*2.3}\approx 0.3$. To show the practicality of each method applied to Lorenz-96, we predict the trajectory of each $X_k$ for a time interval of three times the error doubling time (0.9) after the end of the training set. The measure of inaccuracy used in this paper is the Mean Squared Prediction Error (MSPE) \cite{MSPE1991}.
\begin{equation}
    MSPE = \frac{1}{K}\sum_{k=1}^{K}\frac{1}{I}\sum_{i\in I}(\hat{X}_k(t_i)-X_k(t_i))^2
\end{equation}
where $X_k$ is derived from (\ref{eq:2.0.1}) and (\ref{eq:2.0.2}) and $\hat{X}_k$ is the approximation to $X_k$ derived from a parametrized model (\ref{eq:4.2.1}) through the 4th order Runge-Kutta method, $I$ is the time interval where we compare the trajectories. In this case, $I = [1000, 1000+h, 1000+2h, ..., 1000+0.9]$ where $h=0.01$. The initial conditions are set as: $\hat{X}_k(1000) = X_k(1000), k=1,...,K$ The MSPEs derived for each method are shown in Table \ref{table:MSPE}.

\begin{table*}
\centering
\begin{tabular}{|c c c|}
    \hline
    Method & MSPE & Average K-L Divergence   \\ [0.5ex] 
    \hline\hline
    Regression & 0.03606 & 0.06729 \\ 
    Compressed sensing (Raw) & 5.9993 & 1.9891 \\
    Compressed sensing~[1] (Setting biases to zero) & 0.03398 & 0.10099 \\
    Compressed sensing~[2] (Setting biases to average bias) & 0.03440 & 0.07351 \\
    Compressed sensing~[3] (Adding noise to average bias) & 0.03387 & 0.05919 \\ [1ex] 
    \hline
\end{tabular}
\caption{Evaluation of each method}
\label{table:MSPE}
\end{table*}

\subsection{Comparison of Probability Density Functions}
\label{S:PDFs}
This section will compare the probability density functions (PDFs) of $\hat{X}_k$ and $X_k$ for $k=1,...,K$ during the time interval $t_{test}$. The PDFs are computed using the kernel density estimation, given in (\ref{eq:KernelDensityEstimation}). PDFs are computed for every $X_k$ and $\hat{X_k}$ using their data in the time interval $[500,2000]$.

The Kullback-Leibler divergence was introduced by Kullback and Leibler (1951) and measures the discrepancy between two PDFs. The Kullback-Leibler divergence is calculated as 
\begin{equation}
    D_{KL}(P||Q)=\sum_{x\in X}P(x)\log(\frac{P(x)}{Q(x)})
    \label{eq:KL_Di}
\end{equation}
where $P$ is the PDF calculated from the data of $X_k$ and $Q$ is the PDF calculated from $\hat{X}_k$. We average the Kullback-Leibler divergence in the PDFs of $\hat{X}_k$ and $X_k$ for all $k$ to compare each model. The average Kullback-Leibler Divergence derived for each method is shown in Table \ref{table:MSPE}.

\subsection{Auto-regressive Models in Residuals}
\label{S:AutoRegressive}
So far, we have investigated the deterministic part in stochastic parameterization. The goal of this subsection is to investigate the stochastic part in stochastic parameterization to further improve the regression and compressed sensing models explored so far. This is done by obtaining $\phi,\sigma$ and $\sigma_e$ from the first-order AR model in (\ref{eq:ARMA}). The method of calculating these values is explained in Appendix A. Throughout this subsection, we will refer to the compressed sensing model with noise added to the average bias as the compressed sensing model.

The values for $\phi,\sigma$ and $\sigma_e$ that we obtain are shown in Table \ref{table:ARMA}.
\begin{table*}
\centering
\begin{tabular}{|c c c c c c|}
    \hline
    Method & $\phi$ & $\sigma$ & $\sigma_e$ & MSPE & Ave. K-L Divergence  \\ [0.5ex] 
    \hline\hline
    Regression & 0.9453 & 0.2265 & 0.6945 & 0.02624 & 0.07692  \\ 
    Compressed sensing & 0.9981 & 0.1710 & 2.775 &  0.02066 & 0.07142 \\ [1ex] 
    \hline
\end{tabular}
\caption{AR evaluation of residuals}
\label{table:ARMA}
\end{table*}
Given the values of $\phi$ and $\sigma$, we can run a numerical simulation of the reduced dynamical system with AR residuals. The equation is shown below:
\begin{equation}
    \dot{X}_{k} = - X_{k-1}(X_{k-2} - X_{k+1}) - X_{k} + F - f_k(\textbf{X}(t))+e_k(t_i)
\end{equation}
where $f_k(\textbf{X}(t))$ represents either the Wilk's parameterization~\eqref{WilksPolynomial} or the CS parameterization~\eqref{eq:FinalCSModel}
and $e_k(t_i)$ is generated using the first order AR model below:
\begin{align}
    e_k(t_i)&=\phi e_k(t_{i-1}) + \sigma z_k(t_i) \nonumber \\
    &=\phi e_k(t_{i-1}) + \sigma_e(1-\phi^2)^{1/2} z_k(t_i)
\end{align}
$z_k(t_i)$ is sampled randomly from $N(0,1)$.

Numerical simulations are done using the 4th order Runge-Kutta method. $e_k(t_0)$ is set to $0$ for all $k$. Similar to section 4.2, we evaluate our model by assessing its trajectory as well as its PDF. The results are summarized in Table \ref{table:ARMA}.

A comparison of Table \ref{table:ARMA} with Table \ref{table:MSPE} shows that using an AR model in compressed sensing and regression has significantly improved the prediction of the trajectory of resolved variables as the MSPE is lower. Compressed sensing shows the lowest MSPE and Kullback-Leibler divergence.

The goal of this section was to create a deterministic and AR model to improve the prediction of resolved variables. Based on the results shown, compressed sensing improves the prediction significantly compared to regression. The coefficients of the deterministic model as well as the values of $\phi,\sigma$ and $\sigma_e$ will be useful in the following section where filtering methods will be used to further improve the prediction of resolved variables.

\section{Nonlinear Filtering Applied to Parameterized Models}
\label{S:NonlinearFiltering}
The second component of this paper deals with nonlinear filtering or data assimilation~(DA). 
The main focus of filtering is to combine computational models with sensor data (observations) to predict the dynamics of large-scale evolving systems. 
Filtering provides a recursive algorithm for estimating a signal or state of a random dynamical system based on noisy measurements. 
The signal that is represented by a Markov process cannot be accessed or observed directly and is to be ``filtered" from the trajectory of the observation process which is statistically related to the signal.
Suppose we make a forecast about the behavior at a future time of a complex system with some uncertainties~(randomness) and there is near-continuous data available from remote sensing instrumentation networks. Then, as new information becomes available through observations, it is natural to ask how to best incorporate this new information into the prediction. 
Data assimilation can be simply defined as the combination of observations with the model output to produce an optimal estimate of dynamical model. 

The standard notation for the discrete DA problems from the time $t_{k}$ to $t_{k+1}$ can be formulated as follows:
\begin{equation}
\red
x^{f}(t_{k+1}) = M [x^{f}(t_{k})] + w_{k}
\end{equation}
where $x$ is the model state, $M$ is the model dynamics and $w$ is the model error. The operator $M$ is represented with matrix in discrete case and differential operator in continuous case. 
{\red In this section, for notational simplicity, $x$ represents the model state $\textbf{X}$ as described in previous sections.    }
Observations at time $t_{k}$ can be defined by:
\begin{equation}
\red
y_{k} = H[x^{t}(t_{k})] + v_{k}
\end{equation} 
where $y$ is the observation and $H$ is the observation operator, {\red which is identity $\mathbb{I}$ in our problem}. $v$ is the white noise with associated covariance matrix describing the observation error. 

Data assimilation methods consist of \textit{analysis} (or correction/filtering) and \textit{forecast} (or prediction) steps. At time $t_{k}$, we have the output state of the background forecast $x_{k}^{f}$, and observation state $y_{k}$.
{\red $x^t$ is the truth generated by equation \ref{eq:2.0.1} and \ref{eq:2.0.2}.}
The (linear or nonlinear) analysis step is based on the outputs of $x_{k}^{f}$ and $y_{k}$ and produces analysis $x^{a}_{k}$. Then, applying the model dynamics on the analysis step generates next forecast output as $x_{k+1}^{f}$. \\
\subsection{Ensemble Kalman Filter}
\label{S:EnsembleKF}
It is well documented that particle filters~\cite{Gordon1993} suffer from the well-known problem of sample degeneracy~\cite{Bengtsson2008,Snyder2008}. In contrast, the ensemble Kalman filter~(EnKF)~\cite{Evensen1994,Evensen2003} can handle some of these difficulties, where the dimensions of states and observations are large and the number of replicates is small.
However, EnKF incorrectly treats the non-Gaussian features of the forecast distribution that arise in certain nonlinear systems. 
Nevertheless, we apply ensemble Kalman filter for data assimilation, because they are easy to implement in large systems.

For the Kalman filter~\cite{Kalman60}, the modelling and observation errors are assumed to be independent with Gaussian/normal probability distributions:
\begin{equation}
w_{k} \sim N(0, Q_{k}) \quad\text{and}\quad v_{k} \sim N(0, R_{k})
\end{equation}
where $Q$ and $R$ are the covariance matrices of the model error and observation error, respectively. 
The Kalman filter is of the form:
\begin{equation}
x_{k}^{a} = x_{k}^{f} + K_{k}(y_{k} - H_{k}x_{k}^{f})
\end{equation}
where $K_{k}$ is {\red known as} the \textit{Kalman gain} and $(y_{k} - H_{k}x_{k}^{f})$ is called the \textit{innovation}. After some manipulations, the Kalman gain can be obtained as:
\begin{equation}
K_{k} = P_{k}^{f} H^{T} (H P_{k}^{f} H^{T} + R)^{-1}
\end{equation}
where $P_{k}^{f}$ is the forecast error covariance matrix.

The idea behind EnKF~\cite{Evensen2003} is the propagation and analysis of ensemble members instead of full-rank covariance matrices. Unlike Kalman filter, the EnKF is applicable for nonlinear dynamical models and computationally efficient for high dimensional systems {\red ~\cite{Bocquet2016,Law2015}}. The forecast step of EnKF can be formulated as follows:
\begin{align}
\label{eq:EnKF_forecast}
\red \hat{x}_{j+1}^{n} &\red = \Psi(x_{j}^{n}) + e_{j}^{n}  \hspace{1cm}  \nonumber \\
\red \hat{m}_{j+1} &\red = \dfrac{1}{N} \sum_{n=1}^{N} \hat{x}_{j+1}^{n} \\
\red \hat{C}_{j+1} &\red = \dfrac{1}{N-1} \sum_{n=1}^{N}(\hat{x}_{j+1}^{(n)} - \hat{m}_{j+1})(\hat{x}_{j+1}^{(n)} - \hat{m}_{j+1})^{T} \nonumber
\end{align}
 {\red where $x$ is the model state that is propagated by $\Psi$ and $e_{j}$ is the auto-regressive modeling of the noise as described before}. $m$ is the model state mean, $C$ is the model covariance matrix {\red and $n=1,...,N$ is the number of ensembles}. Then, the analysis step can be performed as follows:
  \begin{align}\label{eq:EnKF_analysis}
S_{j+1} &= H\hat{C}_{j+1}H^{T} + \Gamma \nonumber \\
K_{j+1} &= \hat{C}_{j+1}H^{T}S^{-1}_{j+1} \\
x_{j+1}^{n} &= (I - K_{j+1}H)\hat{x}_{j+1}^{(n)} + K_{j+1}y_{j+1}^{(n)} \hspace{0.5cm}   \nonumber \\
y_{j+1}^{(n)} &= y_{j+1} + \eta_{j+1}^{n}, \hspace{1cm} \nonumber
\end{align}
{\red where $\Gamma$ and $\eta$ are i.i.d's draws from normal distribution and $y_{j}$'s are the observation states. In our problem, the truth is generated by the slow variables of non-parametric Lorenz-96 system, which is given by \eqref{eq:2.0.1} and \eqref{eq:2.0.2}. Then, observations (direct measurements) are produced by adding small sensor noise to the slow variables $X$. And, model used for the propagation of ensemble Kalman filter is either given by compressed sensing or polynomial regression. }
\subsection{Application to Lorenz-96}
\label{S:Autoregressive}
In this section, we compare the filtering results of Wilk's (polynomial) parameterization and compressed sensed parameterization. For each filtering implementation, 40 ensembles are used and observations are included in every 20 time steps. We observe that compressed sensed dynamics provides better prediction results compared to other parameterization methods. There are several softwares that can be used for the filtering problems. Here, we benefit from the open source provided in \cite{PyDA2020} regarding our parameterized Lorenz-96 models. Observations are generated by adding normal noise to the non-parameterized Lorenz-96 model (truth) and filtering trajectories are obtained by using the set of equations given in \eqref{eq:EnKF_forecast} and \eqref{eq:EnKF_analysis} for all cases.
\subsubsection{Autoregressive Wilk's parameterization}
Here, we analyze filtering results of Wilk's parameterization with auto-regressive models. Instead of using random noise for the prediction, we use the auto-regressive modelling of noise as obtained in Section 4.
While keeping the parametrized model same, we include modelled noise to the ensemble dynamics as the following:
\begin{align}
\dot{X}_{k} &=
\G_{k}(X)
- 0.03779 - 0.47362 X_{k} + 0.005197 X_{k}^{2}  \nonumber \\ 
& + 0.004882 X_{k}^{3} - 0.0003142 X_{k}^{4}  + e_{k}(t) 
\end{align}
where $\G_{k}(X)$ is defined in~\eqref{eq:prioriknown} and
\begin{equation}
    e_k(t_i)= 0.9453 e_k(t_{i-1}) + 0.2265  z_k(t_i)
\end{equation}
$z_k(t_i)$ is sampled randomly from $N(0,1)$.
\begin{figure}
\includegraphics[width=0.45\textwidth]{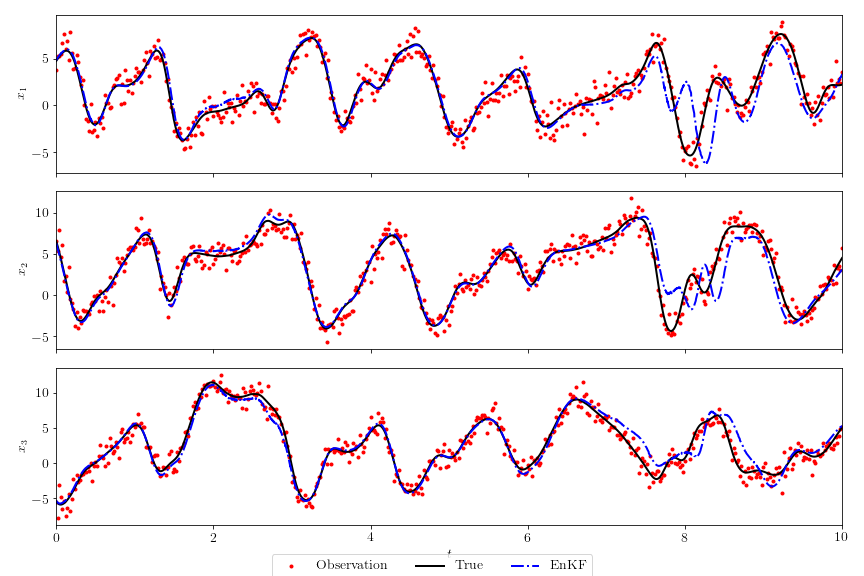}
    \caption{True trajectories, observations and EnKF with Wilk's parametrized model.}
\label{fig:AR_EnKF_Wilks_1}
\end{figure}

Figure \ref{fig:AR_EnKF_Wilks_1} shows the trajectories of truth and observations generated from the truth for the first 3 components of Lorenz-96 model. 
Figure \ref{fig:AR_EnKF_Wilks_2} depicts the truth and prediction (filtering) results of 40 components and MSPE is 0.00940 which is slightly better than deterministic parameterization. 

\subsubsection{Autoregressive Compressed Sensing}
\begin{center}
\begin{figure}
\includegraphics[width=0.5\textwidth]{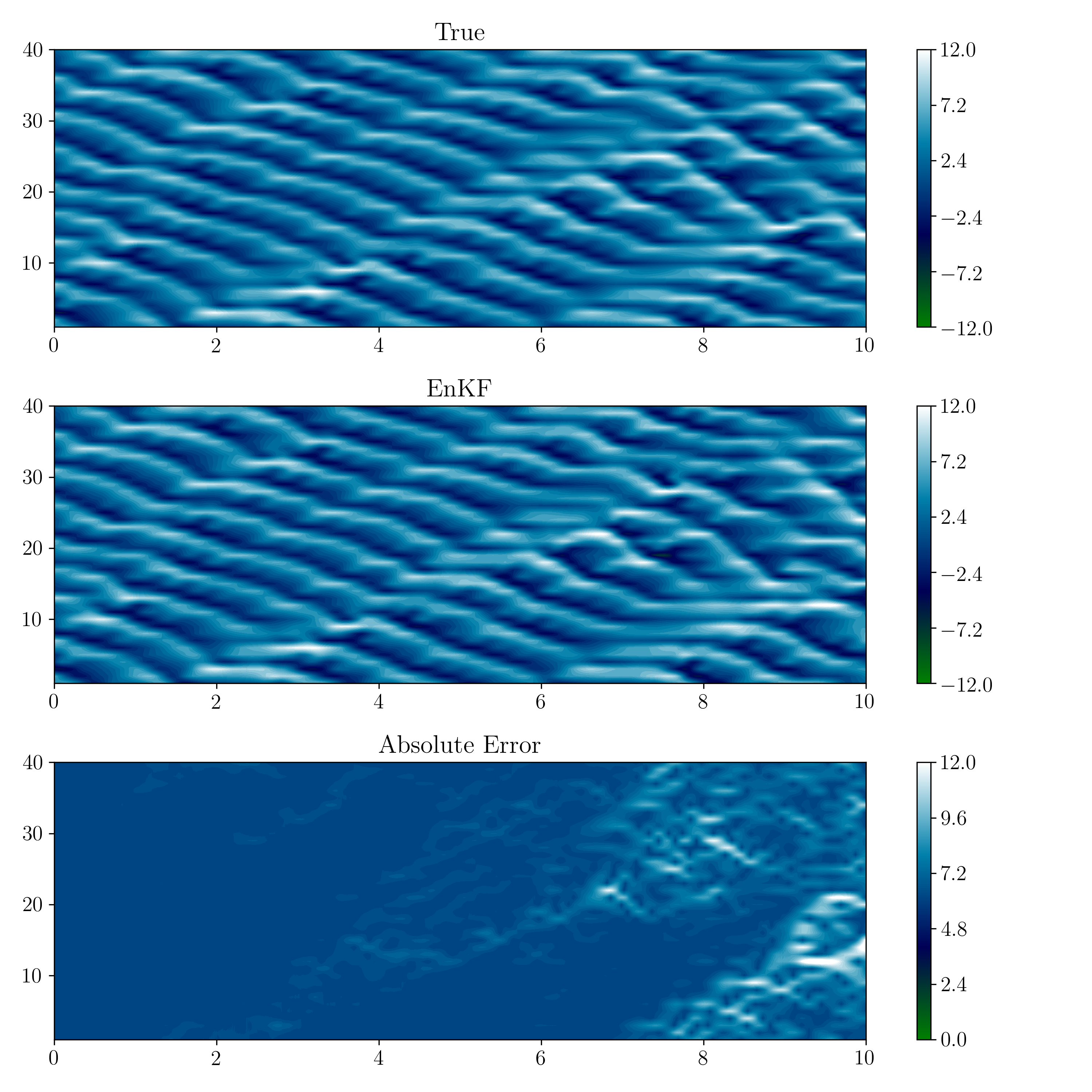}
    \caption{Depiction of Lorenz-96 true and EnKF prediction trajectories for  40 components (Auto-regressive compressed sensing).}
    \label{fig:AR_EnKF_CS_2}
\end{figure}
\end{center}

In addition to compressed sensing method in Section 4, we consider auto-regressive modelling of noise. While keeping the parametrized model same, we include modelled noise to the ensemble dynamics as the following : \\
\begin{align}
\dot{X}_{k} &= 
\G_{k}(X)
- 0.0241 - 0.474 X_{k} + 0.00277 X_{k}^{2} + 0.00509 X_{k}^{3} \nonumber \\ 
& - 0.000306 X_{k}^{4}  + e_{k}(t)
\end{align}
where $\G_{k}(X)$ is defined in~\eqref{eq:prioriknown} and
\begin{equation}
e_k(t_i)= 0.9981 e_k(t_{i-1}) + 0.1710  z_k(t_i)
\end{equation}
$z_k(t_i)$ is sampled randomly from $N(0,1)$.
\begin{figure}
\includegraphics[width=0.45\textwidth]{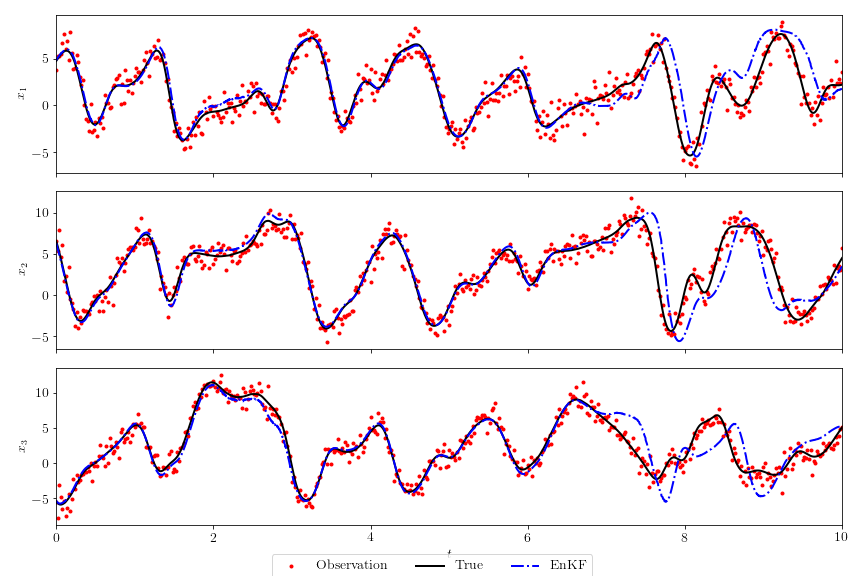}
\caption{True trajectories, observations and EnKF with compressed sensed model.}
\label{fig:AR_EnKF_CS_1}
\end{figure}
Figure \ref{fig:AR_EnKF_CS_1} shows the trajectories of truth and observations generated from the truth for the first 3 components of Lorenz-96 model. 
Figure \ref{fig:AR_EnKF_CS_2} depicts the truth and prediction (filtering) results of 40 components and MSPE is 0.00940 which is slightly better than deterministic parameterization.

\begin{center}
\begin{figure}
\includegraphics[width=0.5\textwidth]{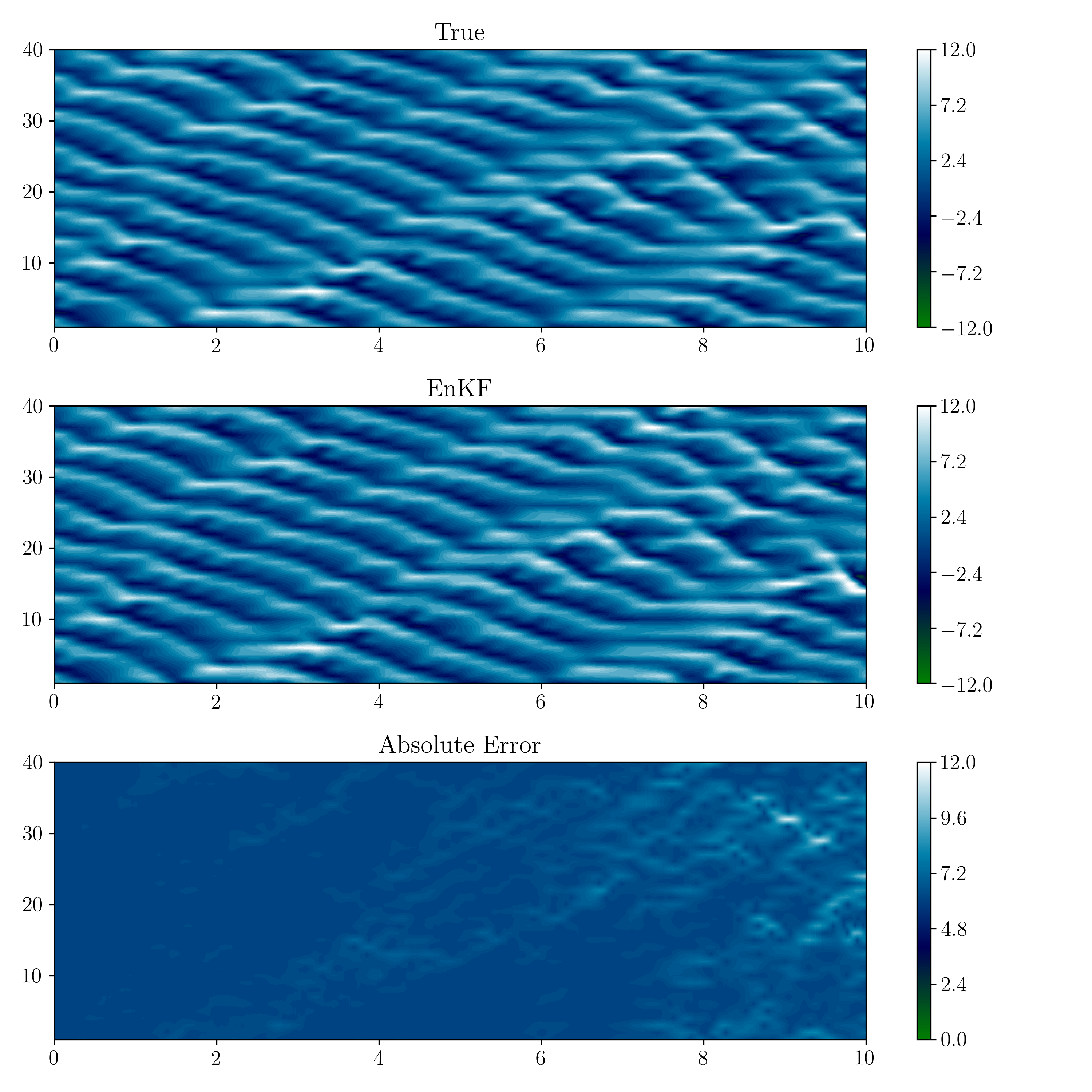}
\caption{Depiction of Lorenz-96 true and EnKF prediction trajectories for  40 components (Auto-regressive Wilk's parameterization) .}
    \label{fig:AR_EnKF_Wilks_2}
\end{figure}
\end{center}

In Figures~\ref{fig:AR_EnKF_Wilks_2} to \ref{fig:AR_EnKF_CS_2}, we see the depictions of filtering results for 40 slow components of Lorenz96 model. 
The first graphs of each figure represent the true trajectories of components where the observations are generated. The second graphs of each figure show the filtering (prediction) results for the system. Finally, the third graphs of each figure depict the L1 errors between the truth and prediction results. All the results about graphs which are shown in this section are discussed in the previous sections related to the individual parameterization methods.
\begin{table}
\label{Table4}
\begin{center}
\begin{tabular}{||c c||} 
\hline
\textbf{Method} & \textbf{MSPE}  \\ [0.5ex] 
\hline\hline
Wilk's parameterization & 0.01006  \\ 
 \hline
 Compressed sensing & 0.01005  \\
 \hline
 Autoregressive Wilk's parameterization & 0.00940 \\
\hline
Autoregressive compressed sensing & 0.00940 \\ [1ex]
\hline
\end{tabular} 
\caption{Parameterization methods and MSPE results including random noise and modelled noise with autoregression}
\end{center}
\label{table:EnKF}
\end{table}

Table 4 shows the parameterization methods with the MSPE errors. As a conclusion, stochastic parameterization methods provide slightly better prediction results comparing to deterministic parameterizations.

\section{Conclusion}
\label{S:Conclusion}
The sheer number of calculations required in directly solving large-scale global weather and climate models becomes computationally overwhelming.
Hence, in this paper, using the multi-scale Lorenz-96 system as a test problem, we introduced a data-driven stochastic-parameterization framework in which the equations of small-scales are integrated and incorporated using compressed sensing to reduce the cost of data assimilation. 
We presented a data-driven, SPR technique, where we learnt to describe the corresponding unknown or unresolved physics from the data due to the fact the governing equations of the unresolved physics are too computationally expensive to solve. 
Since the atmospheric model considered in this paper is inherently chaotic, it has several positive Lyapunov exponents and any small error in the initial conditions grow in finite time.

Random matrices with independent and identically distributed (i.i.d) entries, satisfy RIP with high probability, and the RIP plays an important role in our calculations.
We provided a characterization of the dynamical system based on the Lyapunov exponents for the $F=10$ case.
In this case, the model is strongly chaotic and decorrelation times also suggested the model to be ergodic.
Since the dictionary $\Theta$ was formed as the collection of matrix of all monomials (stored column-wise) of resolved modes and due to the ergodic nature of the simulated time series, the entries satisfied the RIP with high probability. The construction of the optimization problem and computational method were detailed in section~\ref{S:StochasticParameterization}.

For parameterization of $U_{k}$, while previous works mostly focused on suggesting pre-defined polynomial or other model form (or structure) and then estimating parameters, this work is the first attempt to use CS for {\em automatically} and simultaneously discovering polynomial or model form (or structure) and estimating parameters for $U_{k}$ as illustrated by the results in section~\ref{S:StochasticParameterization} and Appendix B.

Specifically, a reduced-order ensemble Kalman filter (EnKF) was applied to the stochastically parametrized model.
For the EnKF, particle locations were corrected based on their respective distances from the observation (innovation). 
Magnitude of the correction was proportional to the error covariance and inverse to the observation noise covariance (Kalman gain). 
In between observations, particles in the EnKF sample evolve according to the original signal dynamics, so particle trajectories are expected to deviate from the truth as time moves further away from the last observation correction. 

Our results show that the compressed sensing model produces the best deterministic parameterization results. After incorporating auto-regression to model residuals, the stochastic parameterizations outperform deterministic parameterizations in terms of prediction accuracy. Compressed sensing with auto-regression produces the best stochastic parameterization results.

Solving for small-scale processes with help of LES models, on a high-resolution grid embedded~(with grid spacing of a few tens of meters and can explicitly simulate the clouds and turbulence represented by parameterization schemes) within the low-resolution grid is very useful, but it's applicability remains limited to its high computational cost.
A long term goal is to create a data-driven scheme, in which LES models embedded in selected grid columns of a global model explicitly simulate subgrid-scale processes which are represented by parameterization schemes in the other columns.
The result is a seamless algorithm which only uses a small scale grid on a few selected columns that simulates interactively within a running global simulation, and is far cheaper than eddy-resolving simulation throughout the global simulation.


\section*{Acknowledgments}
The authors acknowledge partial support for this work from Natural Sciences and Engineering Research Council~(NSERC) Discovery grant 50503-10802, TECSIS /Fields-CQAM Laboratory for Inference and Prediction and NSERC-CRD grant 543433-19. The authors would like to thank Professor Peter Baxendale, Department of Mathematics, University of Southern California, for stimulating discussions, contributions and suggestions on many topics of this paper.


\pagebreak

\section*{Appendix A: Auto-Regression}
Autoregressive (AR) processes are used to predict residuals given a time series of previous residual values.

The first-order AR model is given by (\ref{eq:ARMA}). This can also be written as:
  \begin{equation}
e_k(t_i)=\sigma z_{k,i} + \phi\sigma z_{k,i-1} + \phi^2\sigma z_{k,i-2}+...
\label{eq:ARMA_Multiple}
\end{equation}

Given the time series $e_k(t)$, $\phi$ is used to approximate $e_k(t_i)$ given $e_k(t_{i-1})$. 
$\phi$ is estimated using a least squares minimization process shown below:

 \begin{equation}
C(\phi)=\sum_{k=1}^{K}\sum_{i=2}^{I}(e_k(t_i)-\phi e_k(t_{i-1}))^2
\end{equation}

In order for the process to be stable, $|\phi|<1$ must be satisfied so that the time series $e_k(t_i)$ shown in (\ref{eq:ARMA_Multiple}) does not reach large values.
In~\eqref{eq:ARMA_Multiple}, $\sigma$ is used to measure the variability in $e_k(t_i)$ given $\phi e_k(t_{i-1})$. From calculating $\phi$, we can estimate $\sigma$ by dividing the cost function by the degrees of freedom. This is shown by the following equation:
  
  \begin{equation}
\sigma^2 = \frac{C(\phi)}{K(I-1)-1}
\end{equation}

We then calculate the variance of the process $\sigma_e^2=Var(e_k(t_i))$ using (\ref{eq:ARMA_Multiple}).
\begin{align}
\sigma_e &=\sigma^2 Var(z_{k,i}) + \phi^2\sigma^2 Var(z_{k,i-1}) \nonumber \\
& + \phi^4\sigma^2 Var(z_{k,i-2})+... \nonumber \\ 
& = \sigma^2 + \phi^2\sigma^2 + \phi^4\sigma^2 +... \nonumber \\ 
& \approx \frac{\sigma^2}{1-\phi^2}
\end{align}

The AR model is used in addition to regression and compressed sensing in order to further improve the prediction of resolved variables. 
The estimates for $\phi,\sigma,\sigma_e$ are all important to forming and evaluating the AR model in regression and compressed sensing presented in section~\ref{S:Autoregressive}.

\pagebreak

\begin{table}
\begin{tabular}{c} 
\textbf{Appendix B: Degree-4 Compressed Sensing Model} \\
    \hline
$U_1 = -0.467216X_1 + 0.003429X_1^2 + 0.004481X_1^3  -0.000262X_1^4  -0.024885$\\ [0.5ex]
$U_2 = -0.47689X_2 + 0.0002X_2^2 + 0.005351X_2^3  -0.000299X_2^4  -0.021728$\\[0.5ex]
$U_3 = -0.478369X_3 + 0.002902X_3^2 + 0.005584X_3^3  -0.000358X_3^4  -0.021464$\\[0.5ex]
$U_4 = -0.472069X_4 + 0.002041X_4^2 + 0.005191X_4^3  -0.000307X_4^4  -0.021703$\\[0.5ex]
$U_5 = -0.481931X_5 + 0.001638X_5^2 + 0.005625X_5^3  -0.000345X_5^4  -0.03192$\\[0.5ex]
$U_6 = -0.486479X_6 + 0.001797X_6^2 + 0.006035X_6^3  -0.000385X_6^4  -0.029337$\\[0.5ex]
$U_7 = -0.46724X_7 + 0.001871X_7^2 + 0.004872X_7^3  -0.000281X_7^4  -0.025702$\\[0.5ex]
$U_8 = -0.487255X_8 + 0.002038X_8^2 + 0.00592X_8^3  -0.000368X_8^4  -0.030243$\\[0.5ex]
$U_9 = -0.466759X_9 + 0.003636X_9^2 + 0.004494X_9^3  -0.000258X_9^4  -0.027036$\\[0.5ex]
$U_{10} = -0.463411X_{10} + 0.002444X_{10}^2 + 0.004413X_{10}^3  -0.00024X_{10}^4  -0.013236$\\[0.5ex]
$U_{11} = -0.478315X_{11} + 0.004852X_{11}^2 + 0.004976X_{11}^3  -0.000318X_{11}^4  -0.016471$\\[0.5ex]
$U_{12} = -0.467606X_{12} + 0.004247X_{12}^2 + 0.004642X_{12}^3  -0.00028X_{12}^4  -0.024363$\\[0.5ex]
$U_{13} = -0.46995X_{13} + 0.002823X_{13}^2 + 0.004957X_{13}^3  -0.000297X_{13}^4  -0.018118$\\[0.5ex]
$U_{14} = -0.470231X_{14} + 0.001179X_{14}^2 + 0.004991X_{14}^3  -0.00028X_{14}^4  -0.019603$\\[0.5ex]
$U_{15} = -0.48344X_{15} + 0.003052X_{15}^2 + 0.005955X_{15}^3  -0.00039X_{15}^4  -0.023006$\\[0.5ex]
$U_{16} = -0.470097X_{16} + 0.003608X_{16}^2 + 0.004726X_{16}^3  -0.000275X_{16}^4  -0.030639$\\[0.5ex]
$U_{17} = -0.473192X_{17} + 0.001691X_{17}^2 + 0.005294X_{17}^3  -0.000313X_{17}^4  -0.029741$\\[0.5ex]
$U_{18} = -0.47504X_{18} + 0.004303X_{18}^2 + 0.005322X_{18}^3  -0.000345X_{18}^4  -0.022838$\\[0.5ex]
$U_{19} = -0.473415X_{19} + 0.002968X_{19}^2 + 0.004922X_{19}^3  -0.000289X_{19}^4  -0.011408$\\[0.5ex]
$U_{20} = -0.474502X_{20} + 0.003374X_{20}^2 + 0.005137X_{20}^3  -0.000306X_{20}^4  -0.032649$\\[0.5ex]
$U_{21} = -0.472923X_{21} + 0.003193X_{21}^2 + 0.005178X_{21}^3  -0.000327X_{21}^4  -0.018312$\\[0.5ex]
$U_{22} = -0.466437X_{22} + 0.001616X_{22}^2 + 0.005185X_{22}^3  -0.000318X_{22}^4  -0.017406$\\[0.5ex]
$U_{23} = -0.474267X_{23} + 0.003114X_{23}^2 + 0.005215X_{23}^3  -0.000323X_{23}^4  -0.020402$\\[0.5ex]
$U_{24} = -0.465558X_{24} + 0.002984X_{24}^2 + 0.004562X_{24}^3  -0.000262X_{24}^4  -0.03494$\\[0.5ex]
$U_{25} = -0.474495X_{25} + 0.002997X_{25}^2 + 0.005042X_{25}^3  -0.000302X_{25}^4  -0.018418$\\[0.5ex]
$U_{26} = -0.468641X_{26} + 0.003332X_{26}^2 + 0.004428X_{26}^3  -0.000249X_{26}^4  -0.025618$\\[0.5ex]
$U_{27} = -0.479017X_{27} + 0.002635X_{27}^2 + 0.005268X_{27}^3  -0.000322X_{27}^4  -0.028429$\\[0.5ex]
$U_{28} = -0.466174X_{28} + 0.002651X_{28}^2 + 0.004866X_{28}^3  -0.000289X_{28}^4  -0.034918$\\[0.5ex]
$U_{29} = -0.476855X_{29} + 0.002663X_{29}^2 + 0.005199X_{29}^3  -0.000311X_{29}^4  -0.022651$\\[0.5ex]
$U_{30} = -0.478937X_{30} + 0.001475X_{30}^2 + 0.00543X_{30}^3  -0.00032X_{30}^4  -0.017885$\\[0.5ex]
$U_{31} = -0.45929X_{31} + 0.004295X_{31}^2 + 0.00434X_{31}^3  -0.00026X_{31}^4  -0.024998$\\[0.5ex]
$U_{32} = -0.477136X_{32} + 0.003377X_{32}^2 + 0.005254X_{32}^3  -0.000323X_{32}^4  -0.029709$\\[0.5ex]
$U_{33} = -0.472517X_{33} + 0.001649X_{33}^2 + 0.004832X_{33}^3  -0.000264X_{33}^4  -0.027171$\\[0.5ex]
$U_{34} = -0.470095X_{34} + 0.004248X_{34}^2 + 0.004854X_{34}^3  -0.000301X_{34}^4  -0.018149$\\[0.5ex]
$U_{35} = -0.476457X_{35} + 0.005394X_{35}^2 + 0.005025X_{35}^3  -0.000326X_{35}^4  -0.022066$\\[0.5ex]
$U_{36} = -0.4798X_{36} + 0.000396X_{36}^2 + 0.005354X_{36}^3  -0.000286X_{36}^4  -0.02986$\\[0.5ex]
$U_{37} = -0.468947X_{37} + 0.00297X_{37}^2 + 0.004916X_{37}^3  -0.000289X_{37}^4  -0.031595$\\[0.5ex]
$U_{38} = -0.476784X_{38} + 0.002019X_{38}^2 + 0.005193X_{38}^3  -0.000302X_{38}^4  -0.032644$\\[0.5ex]
$U_{39} = -0.480576X_{39} + 0.004569X_{39}^2 + 0.005093X_{39}^3  -0.000322X_{39}^4  -0.017789$\\[0.5ex]
$U_{40} = -0.482054X_{40} + 0.001322X_{40}^2 + 0.005556X_{40}^3  -0.000332X_{40}^4  -0.020528$\\
    \hline
\end{tabular}
\end{table}

\end{document}